\documentclass[10pt]{article}
\usepackage{amsmath}
\usepackage{amscd}
\usepackage{amsfonts}
\usepackage{amssymb}
\usepackage{mathrsfs}
\usepackage{latexsym}
\usepackage{graphicx}
\usepackage{times}
\usepackage{setspace}
\usepackage{fullpage}
\newcommand{\mbs}[1]{\boldsymbol{#1}}
\title{Semiparametric logistic regression with unknown sizes, 
and its application to bioassays}
\author{Wei Zhang}

\date{}

\DeclareMathOperator*{\argmax}{argmax}
\DeclareMathOperator*{\argmin}{argmin}

\begin{document}
\maketitle

Department of Statistics, 
University of California, Riverside, CA, 92521\\
{wxz118@yahoo.com}
\begin{abstract}
Logistic regression with unknown sizes has many important 
applications in biological and medical sciences. 
All models about this problem in the literature are parametric ones. 
A semiparametric regression model is proposed. 
This model incorporates overdispersion due to  
the variation of sizes, and allows general dose-response 
relations.
An Expectation Conditional Maximization algorithm 
is provided to maximize the log likelihood. The bootstrap 
method can be used to construct confidence intervals for 
regression coefficients. Simulation 
is performed to study the behavior of the proposed model. 
Two real examples are investigated by the proposed model. 
\end{abstract}


\section{Introduction}

Consider that there are $r$ 
observations $(y_i,\boldsymbol{x}_i)$, 
$i=1$, 2, $\dots$, $r$, where $y_i$ is a binomial random variable
 with size $n_i$ and probability $p_i$ and 
$\mbs{x}_i$ is a vector of covariates of length $\varrho$. 
The issue of interest is to investigate how the covariates 
$\mbs{x}_i$ affect the probabilities $p_i$. 
A logistic regression problem arises when the sizes $n_i$ 
are known (e.g., McCullagh and Nelder 1999). It can happen that 
 the sizes $n_i$ are unknown.

The author was motivated to study the logistic regression problem 
with unknown sizes by colony formation assays. These assays 
 are used to assess the cytotoxic effects 
of chemical or physical agents on proliferating cells. 
 In these experiments,  cells are exposed to the agent of interest, 
and then placed onto culture plates for colony formation. 
After some time, visible colonies on each plate are counted to 
decide how many cells survive. 
The initial number of cells put onto each plate 
is usually unknown. Table \ref{table:bovisdata} presents an example 
in which the survival of \textit{M. bovis} cells was studied
 (Trajstman 1989). 
Note that $y_i$ is the number of colonies, and $n_i$ is 
the unknown total number of cells on a culture plate. 

There are many other applications. For example, 
Margolin et al. (1981) studied the effects of 
quinoline on the
number of revertant colonies of Salmonella strain TA98. 
Bailer and Piegorsch (2000) 
reviewed the statistical methods on aquatic toxicology studies
 and took the effect of nitrofen on
the offspring of \textit{C. dubia} as an example. 
 Morton (1981) presented an example of wheat disinfestation by 
hot air. 
Elder (1996) investigated the survival of 
 V79-473 cells and their exposure times to high temperature. 
The radiation damage on jejunal crypts has been studied extensively
(e.g., Khan et al. 1997, Kinashi et al. 1997, Mason et al. 1999, 
Salin et al.  2001, and Goel et al. 2003).

\begin{table}
\small
\caption{The \textit{M. bovis} cell survival data.}
\label{table:bovisdata}
\begin{center}
\begin{tabular}{lccccccccccr}
\hline
 $\%$weight/volume & \multicolumn{10}{c}{No. of \textit{M. bovis} colonies at 
stationarity}
 &sample mean  \\
\hline
 &  &  &  & &  &  &  & &  & &  \\
&\multicolumn{10}{c}{control experiment (no decontaminant)}&\\
  & 52 & 80 & 55 & 50 & 58 & 50 & 43 & 50 & 53 & 54 & 51.8  \\
  & 44 & 51 & 34 & 37 & 46 & 56 & 64 & 51 & 67 & 40 &   \\
  &  &  &  & &  &  &  & &  & &  \\
$\textit{[HPC]}$&\multicolumn{10}{c}{decontaminant: HPC}&\\
 0.75 & 2 & 4 & 8 & 9 & 10 & 1 & 0 & 5 & 14 & 7 & 6.0  \\
 0.375 & 11 & 12 & 13 & 12 & 11 & 13 & 17 & 16 & 21 & 2 & 12.8 \\
 0.1875 & 16 & 6 & 20 & 23 & 23 & 39 & 18 & 23 & 33 & 21 & 22.2  \\
 0.09375 & 33 & 46 & 42 & 18 & 35 & 20 & 19 & 29 & 41 & 36 & 31.9  \\
 0.075 & 30 & 30 & 27 & 53 & 51 & 39 & 31 & 36 & 38 & 22 & 35.7  \\
 0.0075 & 53 & 62 & 38 & 54 & 54 & 38 & 46 & 58 & 54 & 57 & 51.4 \\
 0.00075 & 3 & 42 & 45 & 49 & 32 & 39 & 40 & 34 & 45 & 51 & 38.0  \\
 &  &  &  & &  &  &  & &  & &  \\
$\textit{[Oxalic\,\, acid]}$
&\multicolumn{10}{c}{decontaminant: Oxalic acid}&\\
 5 & 14 & 15 & 6 & 13 & 4 & 1 & 9 & 6 & 12 & 13 & 9.3  \\
 0.5 & 27 & 33 & 31 & 30 & 26 & 41 & 33 & 40 & 31 & 20 & 31.2 \\
 0.05 & 33 & 26 & 32 & 24 & 30 & 52 & 28 & 28 & 26 & 22 & 30.1 \\
 0.005 & 36 & 54 & 31 & 37 & 50 & 73 & 44 & 50 & 37 &  & 45.8  \\
\hline
\end{tabular}
\end{center}
\end{table}

In the literature, the response $y_i$ is usually assumed to be a 
Poisson random variable, such as Wadley (1949) and Margolin et al. (1981). 
 Such an approximation is inappropriate
when $n_i$ and $p_i$ are moderate in size (e.g., Elder et al. 1999). 
Anscombe (1949) considered overdispersion relative to 
the Poisson distribution and developed a model based on the 
negative-binomial distribution. 
Baker et al. (1980) treated $y_i$ as 
a Poisson random variable.  
The $y_i$ in the control group  have a common mean $m$ and  
those in the treatment 
group $m p_i$, where a probit dose-response relation is assumed. 
Trajstman (1989) modified the method of Baker et al. (1980) to 
allow a logistic dose-response 
relation and incorporated overdispersion by assuming a scaled Poisson 
variance-mean relationship.  Morgan and Smith (1992) also based 
their work on Baker et al. (1980), and used a negative-binomial 
variance/mean relationship with a heterogeneity factor to handle 
extra Poisson variation. 
  Kim and Taylor (1994) and Elder et al. (1999)
developed a quasi-likelihood approach by regarding
$y_i|n_i$ as a binomial random variable.
Kim and Taylor (1994) assumed that $E(n_i)=\lambda_i$
and $\text{var}\,(n_i)=\lambda_i\nu$ with $\lambda_i$ 
known and $\nu\geqslant 1$
unknown. Elder et al. (1999) estimated $\lambda=E(n_i)$
with $\text{var}\,(n_i)=\lambda(1+\nu \lambda)$ and $\nu \geqslant 0$.
All previous methods used parametric models.

We propose a semiparametric regression 
 model, in which 
each $n_i$ is assumed to be a Poisson random variable with 
mean $\lambda_i$, and the $\lambda_i$ are assumed to arise 
as a random sample from 
an unspecified mixing distribution. By doing this, 
a rich pool of distributions can be used for $\lambda_i$. 

In Section $2$, a semiparametric model is formulated, and 
an Expectation Conditional Maximization (ECM) algorithm that 
maximizes the log likelihood is described.
  The issues of selecting the number of 
support points and using the bootstrap method are also discussed. 
 Simulation results are shown in Section $3$. 
Section $4$ applies the proposed model to two real examples. 
One is from an \textit{M. bovis} cell survival assay,
and the other from a jejunal crypt stem cell survival assay.
\section{Methods}

\subsection{A semiparametric model}
The probability $p_i$ can be written as
$p_i=h(\mbs{x}_i;  \mbs{\beta})$, where 
$h$ is the inverse of a link function, e.g., 
$h^{-1}=$ logit or probit.  
Note that $h$ is a general function of $\mbs{x}_i$ and $\mbs{\beta}$. 
The unknown size $n_i$ is assumed to be a Poisson random variable with 
mean $\lambda_i$. It is clear that $y_i$ given $\lambda_i$ is a Poisson 
random variable with mean $\lambda_ih(\mbs{x}_i; \mbs{\beta})$. 
The nuisance parameters $\lambda_i$ are further assumed to follow a 
mixing distribution $G$. Because the parameter of interest 
$\mbs{\beta}$ is in the $\varrho$-dimensional Euclidean space, a semiparametric 
regression model arises when $G$ is treated nonparametrically.
The density of a single generic observation $(y,\mbs{x})$ is
\[
f(y;\mbs{x},\mbs{\beta},G)=\int f(y;\mbs{x},\mbs{\beta},
\lambda)dG(\lambda),
\]
where $f(y;\mbs{x},\mbs{\beta},\lambda)$ is a Poisson density 
with mean $\lambda h(\mbs{x}; \mbs{\beta})$, i.e.,
\[
f(y;\mbs{x},\mbs{\beta},\lambda)=\exp\{-\lambda h(\mbs{x}; \mbs{\beta})\}
\{\lambda h(\mbs{x}; \mbs{\beta})\}^y/y!, \quad y=0,1,\dots.
\]
The log likelihood can be written as
\begin{equation}
\label{eq:loglike}
\ell(\mbs{\beta},G)=\sum_{i=1}^{r}\log f(y_i;\mbs{x}_i,\mbs{\beta},G).
\end{equation}

\subsection{An ECM algorithm}

In order to maximize $\ell(\mbs{\beta},G)$ in (\ref{eq:loglike}), 
first we will consider the case that
 $G$ is a discrete distribution with a fixed number of support points. 
Let $G=\sum_{j=1}^{K}\alpha_j\delta(\lambda_j)$, 
where 
$\sum_{j=1}^{K}\alpha_j=1$, $\alpha_j\geq 0$, 
$\delta$ is the indicator function, and $\lambda_j\in (0,\infty)$.
Let $\mbs{\alpha}=(\alpha_1,\, $$\alpha_2,\,\dots,\,\alpha_K)'$,
$\mbs{\lambda}=(\lambda_1,\,$$\lambda_2,\,\dots,\,\lambda_K)'$ and 
$\mbs{\theta}=$($\mbs{\beta},$ 
$\mbs{\alpha}$, $\mbs{\lambda}$). The log likelihood 
$\ell(\mbs{\beta},G)$ in (\ref{eq:loglike}) can be written as
\begin{equation}
\label{eq:disloglike}
\ell(\mbs{\theta})=\sum_{i=1}^r\log \left\{
\sum_{j=1}^{K}\alpha_jf_j(y_i;\mbs{x}_i,\mbs{\beta},
\lambda_j)\right\}.
\end{equation}
One may consider using an EM algorithm to maximize 
$\ell(\mbs{\theta})$ in (\ref{eq:disloglike}).
 However, the M-step in the EM algorithm may be computationally  
unreliable. 

We will consider an ECM algorithm (Meng and Rubin 1993; 
McLachlan and Peel 2000, p148).
The ECM algorithm simplifies the M-step by 
replacing the complicated M-step
 with three computationally simpler and stabler conditional 
maximization (CM) steps.  
It also drives up the log likelihood at each iteration 
 (Meng and Rubin 1993).

Suppose the missing datum is $\mbs{z}=(z_{1},z_{2},...,z_{K})'$, 
the indicator vector for the pair $(\mbs{x},y)$, where 
$z_{j}=1$ for some $j$ and $z_{k}=0$ for all $k\neq j$, i.e.,
 $\lambda=\lambda_j$, $j=1,2,\dots,K$.
Note that $\mbs{z}$ is multinomial distributed with size one and 
probability $\mbs{\alpha}$. The complete density for a single datum 
$(\mbs{x},\mbs{z},y)$ is $
\Pi_{j=1}^{K}\left [\alpha_jf_j(y;\mbs{x},\mbs{\beta},\lambda_j)
\right ]^{z_j}$.
The joint complete log likelihood is
\[
\ell_c(\mbs{\theta})=\sum_{i=1}^r\sum_{j=1}^Kz_{ij}
\big\{\log\alpha_j+\log[f_j(y_i;\mbs{x}_i,\mbs{\beta},\lambda_j)]\big\}.
\]
The expected conditional log likelihood to be maximized
is 
\[
W(\mbs{\theta};\mbs{\theta}^{(0)})=
E_{\mbs{\theta}^{(0)}}\left\{\ell_c(\mbs{\theta})|y_1,y_2,...,y_r\right\}.
\]

The E-step involves getting the conditional expectation of $z_{ij}$, i.e., 
\[
\pi^{(0)}_{ij}=E_{\mbs{\theta}^{(0)}}\big (z_{ij}|y_1,y_2,...,y_r \big)
=\frac{\alpha^{(0)}_jf_j(y_i;\mbs{x}_i,\mbs{\beta}^{(0)},\lambda^{(0)}_j)}
{\sum_{h=1}^K
\alpha^{(0)}_hf_h(y_i;\mbs{x}_i,\mbs{\beta}^{(0)},\lambda^{(0)}_h)}
\]
for 
$i=1,2,\dots,r$ and $j=1,2,\dots,K$.

In the CM-step, we need to maximize the expected conditional complete 
log likelihood
\begin{align*}
W(\mbs{\theta};\mbs{\theta}^{(0)})
&=\sum_{i=1}^r\sum_{j=1}^K 
\pi_{ij}^{(0)}
\log \alpha_j+\sum_{i=1}^r\sum_{j=1}^K \pi_{ij}^{(0)}
\log f_j(y_i;\mbs{x}_i,\mbs{\beta},\lambda_j)\\
&=\text{constant}+\underbrace{\sum_{i=1}^r\sum_{j=1}^K 
\pi_{ij}^{(0)}
\log \alpha_j}_{T_1(\mbs{\alpha)}}\\
&+\underbrace{\sum_{i=1}^r\sum_{j=1}^K\pi^{(0)}_{ij}
\left\{y_i\log\lambda_j+y_i\log h(\mbs{x}_i; 
\mbs{\beta})-\lambda_j h(\mbs{x}_i; 
\mbs{\beta})\right\}}_{T_2(\mbs{\beta},\,\mbs{\lambda})}
\end{align*}
over $\mbs{\alpha},\mbs{\lambda},\mbs{\beta}$ sequentially.
The maximum likelihood estimator (MLE) for $\mbs{\alpha}$ is
 \begin{equation}
\label{eq:alpha}
\alpha_j^{(1)}=r^{-1}\sum_{i=1}^r\pi^{(0)}_{ij}, \,j=1,2,...,K.
\end{equation}
The conditional MLE for $\mbs{\lambda}$
given $\mbs{\beta}=\mbs{\beta}^{(0)}$
is
\begin{equation}
\label{eq:lambda}
\lambda_j^{(1)}=\frac{\sum_{i=1}^r\pi^{(0)}_{ij}y_i}{\sum_{i=1}^r\pi^{(0)}_{ij}
h(\mbs{x}_i; 
\mbs{\beta}^{(0)})},\,
j=1,2,...,K.
\end{equation}
The conditional MLE for $\mbs{\beta}$
given $\mbs{\lambda}=\mbs{\lambda}^{(1)}$
is
\begin{equation}
\label{eq:beta}
\mbs{\beta}^{(1)}=\argmax_{\mbs{\beta}\in \mathcal{R}^\varrho} 
T_2(\mbs{\beta},\mbs{\lambda}^{(1)}).
\end{equation}
Since there 
is no analytic solution for $\mbs{\beta}^{(1)}$ in the optimization 
problem defined in (\ref{eq:beta}), 
a Newton Raphson algorithm is applied. 
The first order derivative of $T_2(\mbs{\beta},\mbs{\lambda}^{(1)})$ is
\[  
\frac{\partial{T_2}}{\partial{\mbs{\beta}}}
=\sum_{i=1}^r\sum_{j=1}^K
\pi_{ij}^{(0)}\left [\frac{y_i}{h(\mbs{x}_i; 
\mbs{\beta})}-\lambda_j^{(1)}\right ]\nabla_{\mbs{\beta}}
h(\mbs{x}_i; 
\mbs{\beta}),
\]
and the second order derivative is
\[
\frac{\partial^2 T_2}{\partial \mbs{\beta} \partial\mbs{\beta}'}
=\sum_{i=1}^r\sum_{j=1}^K \pi_{ij}^{(0)}
\left\{
\left [
\frac{y_i}{h(\mbs{x}_i; \mbs{\beta})}
-\lambda_j^{(1)}
\right ]
\frac{ \partial^2 h(\mbs{x}_i; \mbs{\beta})}
{\partial \mbs{\beta}\partial\mbs{\beta}'}-
\frac{y_i\nabla_{\mbs{\beta}} h(\mbs{x}_i; 
\mbs{\beta})\nabla'_{\mbs{\beta}} h(\mbs{x}_i; 
\mbs{\beta})}{h(\mbs{x}_i;  
\mbs{\beta})^2}\right\}.
\]
The Newton Raphson algorithm is defined by, with
 $\mbs{\beta}_{(0)}=$$\mbs{\beta}^{(0)}$,
\begin{equation}
\label{eq:NRbeta}
\mbs{\beta}_{(t+1)}=\mbs{\beta}_{(t)}-
\left [ \frac{\partial^2{T_2}}
{\partial{\mbs{\beta}}\partial{\mbs{\beta}}'}\Big{|}_{\mbs{\beta}=\mbs{\beta}_{(t)}}
\right ]^{-1}
\left[\frac{\partial{T_2}}
{\partial{\mbs{\beta}}}\Big{|}_{\mbs{\beta}=\mbs{\beta}_{(t)}}
\right].
\end{equation}

\subsection{Selecting the number of support points}

By increasing the number of support points of $G$, the maximized
 log likelihood $\ell(\hat{\mbs{\theta}})$ can be increased. 
One may consider using the global maximizer by trying different values of 
$K$. In order to obtain 
a reasonable and parsimonious fit to the data, 
we propose to choose the number of support points by 
minimizing the BIC (e.g., Wang et al. 1996), i.e., 
\[
\widehat{K}=\argmin_{K\in \{1,2,\dots\}} 
\{-2\ell (\hat{\mbs{\theta}})+\log(r)(2K-1+\varrho)\}.
\]

\subsection{The bootstrap method}
The bootstrap method can be applied to obtain
confidence intervals for the regression coefficients
$\mbs{\beta}$. For a random design, the nonparametric bootstrap 
method can be applied, in which one can sample the pairs $(y_i,\mbs{x}_i)$.
For a fixed design, we propose to use a parametric bootstrap method.
A resample of size $r$ is generated as follows,
\[
y_i^*\sim f(y;\mbs{x}_i,\widehat{\mbs{\beta}},\lambda_i), i=1, 2, \dots, r,
\]
where $ \lambda_i$ is a random variable drawn from
the estimated mixing distribution $\widehat{G}$,
\[
\widehat{G}=\sum_{j=1}^{\widehat{K}}\hat{\alpha}_j\delta(\hat{\lambda}_j).
\]

\section{Simulation}

We report a simulation study in which there is a single covariate 
$x$. There are $10$ replications for each integer $x$ in $[-5,5]$, so that 
$r=110$. A logistic dose-response relation is assumed, i.e., 
\[
\log\left\{\frac{p_i}{1-p_i}\right\}=\beta_0+\beta_1\,x_i.
\]
The intercept $\beta_0$ is fixed to be one.
 A $2^3$ design is considered, i.e.,
\begin{equation*}
\underbrace{\{1,2\}}_{\beta_1}
\times \underbrace{\{(0.5,0.5),(0.25,0.75)\}}_{(\alpha_1,\alpha_2)}
\times \underbrace{\{(100,300),(450,650)\}}_{(\lambda_1,\lambda_2)}.
\end{equation*}

For each setting, $800$ samples are generated.
 The results are shown in Table
\ref{table:simres}. One can observe that the bias, standard deviation
 and mean square error of the slope $\beta_1$ are quite small. The 
 $\beta_1$ falls into the $95\%$ quantile interval, with 
ends being $2.5\%$ and $97.5\%$ quantiles.  

\begin{table}
\caption{Simulation results: sd stands for standard deviation, 
qi for $95\%$ quantile interval, and mse for mean square error.}
\label{table:simres}
\begin{center}
\begin{tabular}{ccccrccc}
\hline
 setting&$\beta_1$ &$(\alpha_1,\alpha_2)$&$(\lambda_1,\lambda_2)$&
bias&sd &qi& mse \\
\hline
1 & 1 &(0.5,0.5) & (100,300)&  0.001 & 0.030 & (0.942, 1.064) & 0.001 \\
2 & 1 & (0.25,0.75)  & (100,300)& 0.003 & 0.025 & (0.954,  1.049) & 0.001 \\
3 & 1 &(0.5,0.5)  & (450,650) & 0.001 & 0.019 & (0.966, 1.040) & 0.000 \\
4 & 1 & (0.25,0.75)  & (450,650) & $-$0.000 & 0.017 & (0.968, 1.032) & 0.000 \\
5 & 2 &(0.5,0.5)  & (100,300) & 0.007 & 0.072 & (1.871, 2.156) & 0.005 \\
6 & 2 & (0.25,0.75)  & (100,300) & 0.007 & 0.063 & (1.901, 2.136) & 0.004 \\
7 & 2 &(0.5,0.5)  & (450,650) & $-$0.000 & 0.045 & (1.919, 2.093) & 0.002 \\
8 & 2 & (0.25,0.75)  & (450,650) & 0.002 & 0.038 & (1.928, 2.076) & 0.001 \\
\hline
\end{tabular}
\end{center}
\end{table}

\section{Example}
\subsection{An \textit{M. bovis} cell survival assay}

The  data in Table \ref{table:bovisdata} are part of Table $1$ in 
Trajstman (1989) and also studied by Morgan and Smith (1992). 
 \textit{M. bovis} cells were treated with one of the decontaminants,  
HPC or oxalic acid with one concentration, then placed on the culture 
plates for colony formation. After $12$ weeks (at stationarity),
 the \textit{M. bovis} colonies were counted.
Trajstman (1989) and Morgan and Smith (1992) treated the count of three
 colonies for HPC dose at $0.00075$ 
 as an extreme observation and omitted it
from all analysis.  However, such a small count can be automatically 
taken care of in the proposed semiparametric model.

An ANOVA model is fitted with a separate factor for each level 
of the decontaminants. 
Let $x_j$ denote a factor for the concentration level $j$ of the  
decontaminants. 
It is assumed that the $p_i$ satisfy that
\begin{equation}
\label{eq:bovisp}
\log\Bigl\{\frac{p_i}{1-p_i}\Bigr\}=\beta_0+\sum_{j=1}^{11}\beta_jx_{ij},
\,\,\, i=1,2,\dots,129,
\end{equation}
where $\beta_0$ is the control effect and $\beta_j$ is the effect
 difference between dose $j$ and the control dose, $j=1,2,\dots,11$.

The results of estimated mixing distributions
 are in Table \ref{table:bovisQ}. 
The smallest BIC corresponds 
to $K=3$. When $K=3$, the estimate $\widehat{G}$ 
 is written as 
\[
\widehat{G}=0.046\,\delta(9.391)
+0.840\,\delta(69.52)+0.115\,\delta(107.1).
\]

\begin{table}
\caption{The estimates of the mixing distribution
and the BIC for the  \textit{M. bovis} data.}
\label{table:bovisQ}
\begin{center}
\begin{tabular}{cccc}
\hline
  component number  & mixing probabilities  & 
support point  &  \\
$(j)$&$(\alpha_j)$&$(\lambda_j)$&BIC\\
\hline
&&$\textbf{one-component mixture}$&\\
 1 & 1 & 71.98 &  1061.1\\
&&$\textbf{two-component mixture}$&\\
 1 & 0.048 & 9.601 & 998.0\\
 2 & 0.952 & 73.59 & \\
&&$\textbf{three-component mixture}$&\\
 1 & 0.046 & 9.391 & \fbox{977.0} \\
 2 & 0.840 & 69.52 &  \\
 3 & 0.115 & 107.1 & \\
&&$\textbf{four-component mixture}$&\\
 1 & 0.045 & 9.376 & 984.2 \\
 2 & 0.180 & 57.53 &  \\
 3 & 0.697 & 73.84 &  \\
 4 & 0.079 & 110.9 & \\
\hline
\end{tabular}
\end{center}
\end{table}

\begin{figure}
\centering
\includegraphics[angle=270,width=12cm]{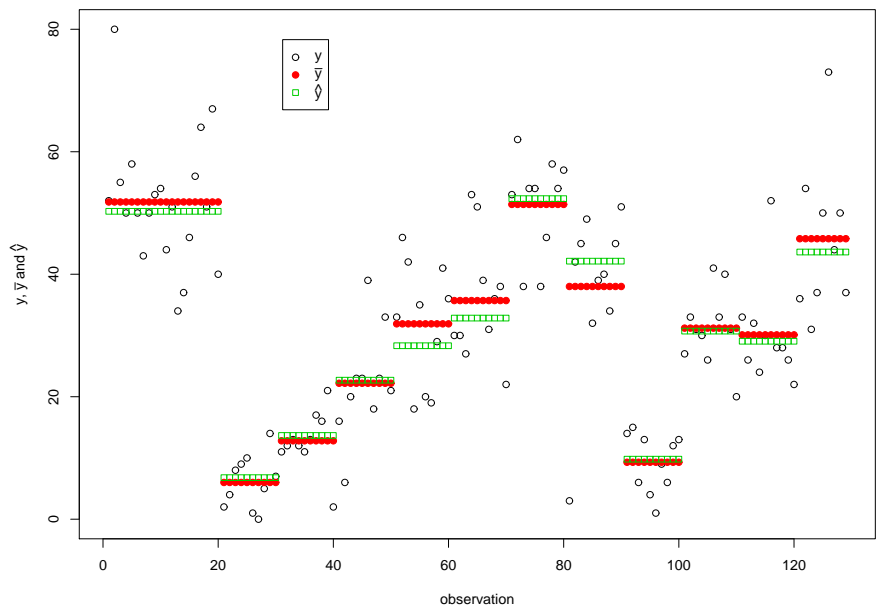}
\caption{The response $y$, sample mean $\bar{y}$ and fitted value $\hat{y}$.}
\label{fig:bovisfit2}
\end{figure}

Table \ref{table:bovisres} presents the results for the regression
 coefficients. In the bootstrap, $200$ resamples are drawn. 
The bootstrap standard errors
 of the regression coefficients are small. Since all $95\%$ confidence 
intervals except those of $\beta_0$ and $\beta_6$ 
 do not include $0$, all treatment 
doses except HPC 0.0075 have  more  negative 
effects on survival of \textit{M. Bovis} cells than the control. 
 The MLEs $\hat{\beta}_6$ and $\hat{\beta}_{9}$
 violate the dose-response monotonicity relationship, i.e., 
increased negative effects on the response
associated with increasing dosage of the 
decontaminants. This is consistent with the 
monotonicity violation in their sample
 means in Table \ref{table:bovisdata}. More investigation is needed for 
the data. The estimates $\hat{\beta}_j$ are not comparable with those 
 in Trajstman (1989) and 
Morgan and Smith (1992), which used a simple linear model in 
(\ref{eq:bovisp}). 
Figure \ref{fig:bovisfit2} presents the responses 
$y$ and their fitted values, which shows that the model fits very well.

\begin{table}
\caption{The estimated regression coefficients,
 bootstrapped standard error, and $95\%$ confidence interval 
for the \textit{M. Bovis} data.}
\label{table:bovisres}
\begin{center}
\begin{tabular}{lcrrc}
\hline
dose& $\beta$&MLE & se & 95$\%$ ci\\
\hline
$\textit{control}$&&&&\\
  &$\beta_0$& 0.882 & 0.117 & (\enspace\, 0.670, \enspace\, 1.125) \\
$\textit{HPC}$&&&&\\
 0.75 &$\beta_1$&  $-$3.131 & 0.209 & ($-$3.615, $-$2.758) \\
0.375 &$\beta_2$&  $-$2.317 & 0.188 & ($-$2.691, $-$1.946) \\
0.1875 &$\beta_3$&  $-$1.639 & 0.180 & ($-$1.983, $-$1.293) \\
0.09375 &$\beta_4$&  $-$1.294 & 0.176 & ($-$1.643, $-$0.960) \\
 0.075 & $\beta_5$& $-$1.034 & 0.193 & ($-$1.443, $-$0.646) \\
0.0075 &$\beta_6$&  0.145 & 0.248 & ($-$0.304,\enspace\, 0.644) \\
 0.00075 &$\beta_7$&  $-$0.506 & 0.196 & ($-$0.857, $-$0.096) \\
$\textit{Oxalic acid}$&&&&\\
5 & $\beta_8$& $-$2.715 & 0.184 & ($-$3.057, $-$2.363) \\
 0.5 &$\beta_9$&  $-$1.155 & 0.191 & ($-$1.533, $-$0.789) \\
 0.05 &$\beta_{10}$&  $-$1.251 & 0.182 & ($-$1.607, $-$0.874) \\
 0.005 &$\beta_{11}$&  $-$0.419 & 0.212 & ($-$0.807,  $-$0.002) \\
\hline
\end{tabular}
\end{center}
\end{table}

\subsection{A jejunal crypt stem cell survival assay}

 Table 1 in Elder et al. (1999) presents a surviving  jejunal crypt 
data set from an experiment done on $126$ mice. Note that the 
colony count of $12$ for dose $9.25$ is redundant and should  
be removed. Kim and Taylor (1994) also investigated the data set.
A jejunal crypt is a compartment containing 
stem cells in a certain 
region of the intestine. These cells are responsible for 
maintaining the function of the intestine. 
In such an experiment, mice are treated by a certain dose of
gamma rays, and then killed to count the number of surviving crypts. 
 Because the experiment needs live mice, the total number of crypts
in each mouse is unknown.
It is assumed that the surviving probabilities $p_i$ satisfy that
\[
\log\left\{\frac{p_i}{1-p_i}\right\}=
\beta_0+\beta_1x_i,\quad i=1,2,\dots,126,
\]
 where $x_i$ is the gamma dose.

The BIC are $724.3$ for $K=1$ and $734.0$ for $K=2$. 
With $\hat{K}=1$, the
 estimated $\widehat{G}$ is degenerated 
at $\hat{\lambda}=196.1$. We draw $200$ bootstrap resamples.
 Table \ref{table:jres} compares the estimates of 
the proposed method with the previous methods. 
All the estimates of previous methods fall into our $95\%$ 
confidence intervals: $(5.089,\, 8.023)$ for $\beta_0$ and
$(-1.241,\, -1.009)$ for $\beta_1$. The standard errors of  
the regression coefficients are 
quite small. Because no confidence intervals include $0$, the 
regression coefficients are significant at the significance level of $0.05$.

\begin{table}
\caption{Jejunal crypt data results from the proposed and
previous approaches 
(logistic regression and Kim's method fix $n_i$ and $E(n_i)$ at 160, 
 respectively; Kim's and Elder's quasi-likelihood method of moments estimates
come from Elder et al. (1999)).}
\label{table:jres}
\begin{center}
\begin{tabular}{c|ccccc}
\hline
& \multicolumn{4}{c}{estimate (standard error)}\\
\cline{2-5}
&logistic &Kim's&
 Elder's&proposed\\
\hline
$\beta_0$ & 7.432 (0.175)&7.410 (0.191) &6.727 (0.725)&6.705 (0.746)
\\
$\beta_1$ & $-$1.185 (0.024) &$-$1.183 (0.026)& $-$1.126 (0.061)&$-$
1.124 (0.059)\\
$\lambda$ &--- &--- & 194.7 (43.4)&196.1\\
\hline
\end{tabular}
\end{center}
\end{table}

\section{Discussion}

We propose a flexible semiparametric model for the logistic  regression 
problem with unknown sizes, in which the regression coefficients 
can be estimated together with the nuisance parameter, the mixing 
distribution.

The parameter estimates in the proposed model can be obtained effectively 
by an ECM algorithm. 
When one runs the ECM algorithm, 
good initial values will help 
find the MLEs  quickly. 
One may run a Poisson regression analysis 
to find the initial values of $\mbs{\beta}$.

\nocite{*}
\bibliography{ref}
\bibliographystyle{apalike}
\end{document}